\documentclass[12pt]{article}
\usepackage{amsmath,amsfonts,latexsym}

\numberwithin{equation}{section}

\newcommand{\C}{{\mathbb C}}

\newcommand{\Z}{{\mathbb Z}}

\newcommand{\G}{{\cal G}}

\newcommand{\n}{{\underline{n}}}
\newcommand{\s}{{\underline{s}}}
\newcommand{\m}{{\underline{m}}}

\begin{document}
\begin{center}
{\LARGE{\bf Irreducible Representations for Toroidal Lie-algebras}} \\
[7cm]
{\bf S. Eswara Rao} \\
{\bf School of Mathematics} \\
{\bf Tata Institute of Fundamental Research} \\
{\bf Homi Bhabha Road} \\
{\bf Mumbai - 400 005}
 {\bf India} \\ [5mm]
{\bf email: senapati@math.tifr.res.in}
\end{center}

\pagebreak

\paragraph*{Introduction}

Let $\G$ be simple finite dimensional Lie algebra over complex numbers $\C$.
The universal central extension of ${\cal G} \otimes \C [t_1^{\pm 1},
\cdots,t_n^{\pm 1}]$ are called toroidal Lie algebras.  For $n=1$ they are
precisely 
the affine Kac-Moody Lie-algebras so that toroidal Lie-algebras are generalisations of
affine Kac-Moody Lie-algebras.  The major difference in the toroidal case is
that the centre is infinite dimensional unlike the affine case where the
centre is one dimensional.  This poses a major difficulty in studying 
representation theory of toroidal Lie algebras.  The toroidal Lie-algebras
are naturally $\Z^n-$graded and there is an infinite dimensional
centre with non-zero degree.  Generally in an irreducible representation the
non-zero degree central operator does not act as scalars but acts as a invertible central operators.  The main purpose of this paper is to prove that the
study of irreducible representation for toroidal Lie-algebra is reduced to the study of irreducible representation where the infinite centre acts as
scalars.

We now give a more details of the paper.  We  fix an integer $ n \geq 2$.  Let
$A_n = \C [t_1^{\pm 1}, \cdots , t_n^{\pm 1}]$ be a Laurent polynomial in
$n$
commuting variables.  We consider the universal central extension of $\G
\otimes A_n$ and add degree derivations $D= \{d_1, \cdots, d_n\}$ which we
call the toroidal Lie-algebra $\tau$.  (see section 1).  Let $V$ be 
irreducible representation for the toroidal Lie-algebra with finite 
dimensional weight spaces.  We prove in section 1 that such a representation
is actually a representation for $\G_{{\rm aff}} \ \otimes A_{n-1} $ up to
an automorphism of $\tau$ (see 1.12).

In the process we prove an important result which is independent interest.
Let $\overline{\tau}$ be the quotient of $\tau$ by non-zero degree central
elements.  Then $\overline{\tau}$ does not admit representations with finite
dimensional weight spaces where some zero degree central operator acts
non-trivially (Corollary (1.11)).  See [BC], [MJT] and [F] where
representations for $\overline{\tau}$ are studied.

We prove in section 2 that the above irreducible representation admits an
irreducible quotient $V_1$ for $\G_{{\rm aff}} \otimes A_{n-1} \oplus
\C d_n$ (Theorem (2.5)).  Note that the quotient is only graded by
$d_n$ which comes from the gradation of $\G_{{\rm aff}}$.  In this
quotient the central operators act by scalars (Lemma (3.1)).

In the last section we give a $\G_{{\rm aff}} \otimes A_{n-1} \oplus D$
module
structure on $V_1 \otimes A_{n-1}$ and prove complete reducibility in
Proposition (3.8).  In Theorem (3.9) we prove that the original module
$V$ is isomorphic to a component of $V_1 \otimes A_{n-1}$.  In fact all
components are isomorphic up to a grade shift.  Thus we conclude that
 the study of irreducible representations for $\tau$ with finite
dimensional weight spaces is reduced to the study of irreducible
representation for $\G_{{\rm aff}} \otimes A_{n-1} \oplus \C d_n$ with
finite dimensional weight spaces where the centre acts as scalars.

It may be appropriate to mention that irreducible integrable representations
 for $\tau$ with finite dimensional weight spaces are classified by the
author in [E3] and [E4].  See also [Y].  The integrable modules are not
completely reducible. Indecomposable modules which are not necessarily
irreducible has been studied by [CT].  Further vertex representations
are constructed and studied in [MEY], [EM], [E1], [JMT] and [CO].

\paragraph*{Section 1}

Let $\G$ be simple finite dimensional Lie algebra over complex numbers
$\C$.  Let (,) be $\G$-invariant non-degenerate symmetric billinear form.
Fix a positive integer $n \geq 2$.    Let $A=A_n =\C[t_1^{\pm 1}, \cdots,
t_n^{\pm 1}] $ be a  Laurent polynomial ring in $n$ commuting variables.
For $\underline{m}=(m_1, \cdots, m_n) \in \Z^n$, let $t^{\underline{m}}=
t_1^{m_1} \cdots t_n^{m_n} \in A$.  For any vector space $V$, let $V_A=
V \otimes A$.  Let $v (\underline{m}) =v \otimes t^{\underline{m}}$.  Let
$Z=\Omega_A / d_A$ be the module of differentials.  That is $Z$ is
spanned by vectors $t^{\underline{m}} K_i, 1 \leq i \leq n, \underline{m}
\in \Z^n$ with a relation.

\paragraph*{(1.1)} $\sum m_i t^{\underline{m}} K_i =0$.  

Let $\tau = \G \otimes A \oplus Z \oplus D$ where $D$ is spanned by $d_1,
\cdots, d_n$.  We
 will now define a Lie algebra structure on $\tau$.

\paragraph*{(1.2)}  $F$ or $X,Y \in \G, \underline{r}, \underline{s} \in \Z^n$.
\begin{enumerate}
\item[(1)] $[X (\underline{r}), Y (\underline{s})] =[X,Y] (\underline{r}+ 
\underline{s})+ (X,Y) d (t^{\underline{r}}) t^{\underline{s}}$ \ where
$d (t^{\underline{r}}) t^{\underline{s}} = \sum r_i t^{\underline{r}+s}K_i$.
\item[(2)] $Z$ is central in $\G \otimes A \oplus Z$.
\item[(3)] $[d_i, X(\underline{r})]=r_i X (\underline{r}),  \ [d_i,
d(t^{\underline{r}}) t^{\underline{s}}]= (r_i+s_i)d (t^{\underline{r}})
t^{\underline{s}}, [d_i, d_j]=0$.
\end{enumerate}

\paragraph*{Proposition (3.1)} ([MY], [K]).  $\G \otimes A \oplus Z$ is
the universal central extension of $\G \otimes A$.
\paragraph*{(1.4)}  Let $h$ be a Cartan subalgebra of $\G$.  Let
$\tilde{h}$
be the span of $h, K_1, \cdots K_n$ and $d_1, \cdots, d_n$.  We call
$\tilde{h}$ a Cartan subalgebra of $\tau$.  Let $\delta_1, \cdots, \delta_n$
be in $\tilde{h}^*$ defined by $\delta_i (h) =0, \delta_i (K_j)=0$ and
$\delta_i (d_j)= \delta_{ij}$.  For $\underline{r} \in \Z^n, \delta_{\underline{r}} = \sum r_i \delta_i$.  Let $\stackrel{\circ}{\triangle}$ be the root 
system of $\G$ and let $\G= \displaystyle{\oplus_{\alpha \in 
\stackrel{\circ}{\triangle}}} \ \G_{\alpha} \oplus h$ be the root space 
decomposition.  Let $\triangle =\{\alpha+\delta_{\underline{r}},
\delta_{\underline{s}}, \ \alpha
\in \stackrel{\circ}{\triangle}, \underline{r}, \underline{s} \in \Z^n \}$.
Let
$$
\begin{array}{lll}
\tau_{\alpha+ \delta_{\underline{r}}} &=& \G_{\alpha} \otimes t^{\underline{r}}, \alpha \in \stackrel{\circ}{\triangle}, \underline{r} \in \Z^n \\
\tau_{ \delta_{\underline{r}}} &=& h \otimes t^{\underline{r}}, \underline{r}
\in \Z^n, \underline{r} \neq 0 \\
\tau_0 &=& \tilde{h}
\end{array}
$$
Then clearly $\tau= \displaystyle{\oplus_{\alpha \in \triangle}} \tau_{\alpha}$
is a  root space decomposition of $\tau$.

\paragraph*{(1.5)} Let $G=GL (n, \Z)$ be the set of $n \times n$
matrices with entries in $\Z$ whose determinant is $\pm 1$. For every $B$ in 
$G$ we define an automorphism of $\tau$ again denoted by $B$.
Let $\underline{r}, \s \in \Z^n$ and let $\underline{m}^T = B \underline{r}^T$
 and $\n^T = B \s^T$.
Define
$$BX (\underline{r}) = X (\underline{m}), X \in \G$$
$$Bd (t^{\underline{r}}) t^{\s} =d  (t^{\underline{m}}) t^{\n}$$
Let $(d^1_1, \cdots, d_n^1)^T =A^{-1} (d_1, \cdots d_n)^T$.  Define
$B d_i =d_i^1$.  It is easy to check that $B$ is an automorphism of
$\tau$.

\paragraph*{(1.6)}  We recall some simple results from [E4].

Let $V$ be irreducible weight module for $\tau$ with finite dimensional
weight spaces with respect to $\tilde{h}$.  the following is proved in
[E4].

\paragraph*{Lemma (1.7)} (Lemma (4.3) and (4.4) of [E4]).  Let $z$ be in
$\Z$ of degree $\m$.

(1) ~ Suppose $z v \neq 0$ for some $v$ in $V$.  Then $zv \neq 0$ for all 
non-zero $v$ in $V$.

(2)~~ Suppose $z \neq 0$.  Then there exists a central operator $T$
on $V$ such that $zT = Tz =Id$.

Let $L= \{\underline{r}, - \underline{r} \mid t^{\underline{r}}K_i \neq 0$ on $V$ for some $i \}$.

Let $T$ be the $\Z$ linear span of $L$.  Clearly $T$ is a subgroup of $\Z^n$.
Let $k$ be the rank of $T$.  Clearly $k$ is invariant under the automorphism
defined above.  Now by standard basis theorem there exists a $\Z$
basis $\m_1, \cdots, \m_n$ of $\Z^n$ such that $k_1 \m_1, \cdots,
k_k \m_k$ is a basis of $T$ for some non-zero integers $k_1, \cdots k_k$.  It is a standard fact that there exists a $B$ in $G$ such that $B \m_i =(0,
\cdots 1, 0 \cdots 0)$  (1 in the $i$th place).  Thus up to automorphism,
we can assume that there exists non-zero central operator $z_1, \cdots
z_k$ on $V$ such that  their degree are 
 $(k_1, 0 \cdots 0) \cdots (0, \cdots k_k, 0 \cdots 0)$.  We have the 
following:

\paragraph*{Proposition 1.8} ~~ (Theorem (4.5) of [E4]).  Let $V$ be 
irreducible module for $\tau$ with finite dimensional weight  spaces with 
respect to $\tilde{h}$.  Then up to automorphism the following  is true on $V$.
\begin{enumerate}
\item[(1)] There exists non-zero central operators $z_1, \cdots, z_k$ of
degree \\
 $(k_1, 0 \cdots 0) \cdots (0 \cdots k_k, 0 \cdots 0)$.
\item[(2)]  $k<n$
\item[(3)] $t^{\underline{r}} K_i =0$ for all $i$ and for all
$\underline{r}$ such that for some $j, k+1 \leq j \leq n, \  r_j \neq 0$.

\item[(4)] $t^{\underline{r}} K_i =0 ~~ 1 \leq i \leq k, \forall 
\underline{r} \in \Z^n$.
\item[(5)] There exists a proper submodule $W$ of $\G \otimes A \oplus Z 
\oplus D_k$ (where $D_k$ is spanned by $d_{k+1}, \cdots d_n)$ such that $V
/W$ has finite
dimensional weight spaces with respect to $h \oplus \sum \C K_i \oplus
D_k$. 
\end{enumerate}

\paragraph*{Theorem (1.9)} Let $V$ be an irreducible $\tau$-module with
finite dimensional weight spaces with respect to $\tilde{h}$. Let $k$ be as 
defined above.   Suppose some central operator is non-zero.  Then
$k=n-1$.  

Before we prove the theorem, we prove an important proposition which is
also of independent interest.

Let $H$ be a finite dimensional vector space over $\C$ with a  non-degenerate 
billinear form $<,>$.  Consider the Lie-algebra $\tilde{H}= H \otimes
\C [t^{\pm 1}_{1}, t_2^{\pm 1} ] \oplus \C C_1 \oplus \C C_2 \oplus \C d_1
\oplus \C d_2$ with the following Lie bracket
\begin{enumerate}
\item[(1)] $[h_1 \otimes t_1^{r_1} t_2^{r_2}, h_2 \otimes t_1^{s_1} t_2^{s_2}]
=(h_1, h_2) \delta_{r_1+s_1,0} \delta_{r_2+s_2,0} (r_1 C_1+r_2 C_2)$.
\item[(2)]$C_1, C_2 $ are central
\item[(3)] $[d_i, h_1 t_1^{r_1} t_2^{r_2} ] = r_i h_1 t_1^{r_1} t_2^{r_2}, \
[d_i, d_j]=0$.
\end{enumerate}

Let $h =H \oplus \C C_1 \oplus  \C C_2 \oplus \C d_1 \oplus \C d_2$. 

\paragraph*{Proposition (1.10)}  Let $W$ be a $h$-weight module for $\tilde{H}$
 with finite dimensional weight spaces.  Suppose, at least on $C_i$ acts
non trivially.  Then such a $W$ does not exists.

\paragraph*{Proof} The proof is inspired by Futorny's work [F].

First note that $GL(2, \Z)$ acts on $\tilde{H}$ as automorphism.  Now twisting
 by automorphism we can assume that both central elements act non
trivially on $W$.

Let $W=\displaystyle{\oplus_{\n \in \Z^2}} \ W_{\n}$ be the weight space
decomposition.  Let $h t_1^{r_1} t_2^{r_2} = h \otimes t_1^{r_1} t_2^{r_2}$.

Choose $h \in H$ such that $(h,h) =1$.  Let $e=h t_1 ht_1^{-1}$ and $f=
h t_2 h t_2^{-1}$.  Then $e$ and $f$ leaves each finite dimensional
 space $W_{\n}$ invariant.

Fix $\n \in \Z^2$ such that
$W_{\n} \neq 0$.  Let $v$ be in $W_{\n}$ such that, $v$ is a common eigen
vector for $e$ and $f$.  Let $v$ be a submodule of $W$ generated by $v$.
Let $V_{\m}= V \cap W_{\m}$.  We call a vector of the form $h_1
t_1^{r_{11}} t_2^{r_{21}} \cdots h_n t_1^{r_{n1}} t_2^{r_{n2}} v$ a 
monomial.  Note that each monomial is a eigen vector for $e$ and $f$.
Conversely any eigen vector of $e$  or $f$ is a linear combination
of monomials with same eigenvalue.  From Futorny's work (Lemma (4.2) of
[F]) it follows that the eigen values of $e$ and $f$  are of the form
$m_1 C_1$ and $m_2 C_2$ for some $m_i \in \Z$.  Choose a eigen
 vector $w$ in $V_{\n} $ of eigenvalue  $m_2 C_2$ for $f$ such that the
absolute value of $m_2$ is maximal.

Suppose $m_2 >0$.  From the above observations we can assume $w$ is a monomial.
Suppose $ew=mw$.

\paragraph*{Claim 1}  Suppose $m>0$.  Then $ht_1^{\ell} t_2^k w=0$ for all 
$k \geq 2$ and for all $\ell > 0$.  Suppose $ht_1^{\ell} t_2^k w \neq 0$
for some $k \geq 2$ and for some $\ell >0$.  Consider $h t_2 h t_2^{-1} h
t_1^{\ell} t_2^k w= m_2 C_2 h t_1^{\ell} t_2^k w$.
We are using the fact that $k \geq 2$ so that $f$ commutes with
$ht^{\ell}_1 t_2^k$.
This proves that $ht_2^{-1} ht_1^{\ell} t_2^k w \neq 0$.  Let  $w_s =
(h t_2^{-1})^s ht_1^{\ell} t_2^k w$.  We prove by induction on
$s$ that $w_s \neq 0$.  Consider  $h t_2 w_{s+1}= fw_s = (s+m_2)w_s$.
By induction $w_s \neq 0$ and we know $s+m_2 >0$.  This proves $w_{s+1}
\neq 0$.  In particular $w_k \neq 0$.

\paragraph*{Sub-claim}  $(h t_1^{-1})^{\ell} w_k \neq 0$.  Consider 
$ht_1 ht_1^{-1} w_k=m w_k$.  Since $m \neq 0$, this proves $ht_1^{-1} w_k 
\neq 0$. Consider $v_s = (h t_1^{-1})^s w_k$.  We prove by induction on
$s$ that $v_s \neq 0$.  Consider $ht_1 v_{s+1} = e v_s = (s+m) v_s$.  By
induction $v_s \neq 0$ and we know that $s+m>0$.  This proves $v_{s+1}
\neq 0$.  In particular $v_{\ell}  = (h t_1^{-1})^{\ell} w_k \neq 0$.
But $f v_{\ell} = (k+m_2) v_{\ell}, v_{\ell} \in V_{\n}$.  This is a
contradiction to the
maximality of $m_2$. This proves the claim 1.

\paragraph*{Claim 2.}  suppose $m \leq 0$.  Then $ht_1^{\ell}
 t_2^k w=0$ for all $k \geq 2$ and for all $\ell < 0$.  Consider $ht_1^{-1}
ht_1 w = ht_1 \cdot ht_1^{-1} w -C_1 w$ 
$$= (m-1) C_1 w$$.

Note that $m-1 <0$.  By arguing as in the claim 1 we prove that
$$(h t_2^{-1})^k ht_1^{\ell} t_2^k w \neq 0.$$  Now we claim that
$$b= (ht_1)^{- \ell}  (h t_2^{-1})^k ht_1^{\ell} t_2^k  w \neq 0.$$  This
can be
established by arguing as in the claim 1.  Now using $ht_1^{-1} ht_1$
instead of $ht_1 ht_1^{-1}$, noting $m-1 <0$.  Now $fb = (m_2+ k)b $.  It
is a contradiction to the maximality to $m_2$.

\paragraph*{Claim 3}  $ht_1^{\ell} t_2^{-k} w \neq 0$ for all $k \geq 2$ 
and for all $\ell > 0$ or for all $\ell <0$ such that $\ell C_1 \pm
k C_2 \neq 0$.  This follows from claim 1 and 2 by applying $ht_1^{- \ell}
t_2^k$ on the above vector.

\paragraph*{Case 1}  Assume $ht_1^{\ell} t_2^{-k} w \neq 0$ for all
 $k \geq 2$ and all $\ell \geq 1$.  Suppose $ht_1^{- \ell} t_2^{-k} w
\neq 0$ for some $k \geq 2$ and for all $\ell > 0$ such that
$\ell C+1 \pm k C_2 \neq 0$.  Then we claim that
$$\{h t_1^{- \ell} t_2^{-k} ht_1^{\ell} t_2^{-k} w, \ell C_1 \pm kC_2 
\neq 0\}$$ is linearly independent.  Suppose the claim is true.  Then 
$V_{\n+ (0, 2k)}$ is infinite  dimension a contradiction.  This prove the 
proposition.

Suppose
$$\sum a_{\ell_i} h t_1^{-\ell_i} t_2^{-k} ht_1^{\ell_i} t_2^{-k} w= 0.$$
Apply
$ht_1^{- \ell_1} t_2^k$ to the above sum.  It commutes with the first
term.  Thus we get
$$a_{\ell_1} ht_1^{- \ell_1} t_2^{-k} (\ell C_1 -k C_2) w.$$
This forces from our assumptions that $a_{\ell_1}=0$.  Similarly $a_{\ell_i }
=0$.  Thus the above set is linearly independent.  Now suppose $ht_1^{-\ell} 
t_2^{-k} w =0$ for all $k \geq 2$  and for some $\ell > 0 \ni \ell C_1 \pm
 kC_2 \neq 0$.  Then it follows that $ht_1^{\ell} t_2^k w \neq 0$.  Now
consider the set
$$\{ht_1^{\ell} t_2^k ht_1^{\ell} t_2^{-k} w, \ell C_1 \pm k C_2 \neq
0\}$$
which can be shown to be linearly independent.

This contradicts the finite dimensionality of the weight space.

\paragraph*{Case 2}   $ht_1^{\ell} t_2^{-k} w \neq 0$ for all $k \geq 2$ and 
for all $\ell < 0$.  This case can be dealt with as in the  Case 1.

We are now left with the case $m_2 \leq 0$.  Consider 
$$ht_2^{-1} ht_2 w=
h t_2 h t_2^{- 1} w -C_2 w$$
$$= (m_2-1)C_2 w.$$
Now we have $m_2 -1 < 0$.  Now the above proof works with $ht_2^{-1} ht_2$ 
instead of $ht_2 ht_2^{-1}$. \hfill{[QED]} \\
We first deduce an interesting corollary.   Let $\overline{\tau} =\G
\otimes A \oplus \sum \C K_i \oplus \sum \C d_i$ be a Lie-algebra defined
by
$$[X (\underline{r}), Y (\s)]=[X,Y] (\underline{r}+\s)+(X,Y) 
\delta_{\underline{r}+\s,0} \sum r_i K_i$$
$K_i$ are central.
$$[d_i, X(\underline{r})]=r_i X (\underline{r}), \ [d_i, d_j]=0.$$

\paragraph*{Corollary (1.11)} Let $V$ be a weight module for
$\overline{\tau}, n \geq 2$ generated by a weight vector with finite
dimensional weight spaces.
Assume at least one $K_i$ acts non trivially on the generator.  Then such
a $V$ does not exist.

\paragraph*{Proof}  Follows from the above proposition by restricting the 
module to the subalgebra generated by $h \otimes t_i^r t_{i \pm 1}^{s}$.

Proof of the Theorem (1.9).

>From Proposition (1.8) we know that $k<n$.  Suppose $k<n-1$.

\paragraph*{Case 1}  Assume $K_i \neq 0$ on $V$ for some $i$.
Then $i \geq k+1$ from Proposition (1.8) (4).  Let us say $i=k+1$.  Then
$k+2 \leq n$.  Consider the Lie sub algebra generated by
$ht^r_k t^s_{k+1}$.  Consider $[h_1 t_k^{r_1} t_{k+1}^{s_1}, h_2 t_k^{r_2}
t_{k+1}^{s_2}]= (h_1, h_2) \delta_{r_1+r_2,0} \delta_{s_1+s_2,0} (r_1
K_{k+1}1+s_1 K_{k+2})$.  This is because $d(t_k^{r_1} t_{k+1}^{s_1})
t_k^{r_2}
 t_{k+1}^{s_2}=0$ unless $r_1+ r_2=0$ and $s_1+s_2=0$ from Proposition
(1.8) (3).  Now by restriction the module to the above sub algebra, we get
a contradiction to Proposition (1.10).

\paragraph*{Case 2}  Assume $t^{\underline{m}} K_i$ is non zero on $V$.
Then by Proposition (1.8) (3), $i \geq k+1$.  Further $m_{k+1} =0 \cdots
=m_n$.  Let us say $i=k+1$.  Consider the subalgebra generated by \\ 
(1) ~ $h t^{\m} t_k^r t_{k+1}^s, s > 0$ or $s=0$ but $r>0$, \\
(2) $ ht_k^r t_{k+1}^s, s<0$ or $s=0, r<0$ \\
(3) $h t^{\m}$.

Consider in the above sub algebra,
$$[h_1 t^{\m} t_k^{r_1} t_{k+1}^{s_1}, h_2 t_k^{r_2} t_{k+1}^{s_2}]$$
$$=(h_1, h_2) \delta_{r_1+s_1,0}  \delta_{r_2+s_2,0} (r_1 t^{\m} K_{k+1}+r_2
t^{\m} K_{k+2}).$$
All other brackets of the above sub algebra are zero. This is due to the 
Proposition (1.8).  It is easy to see that the above sub algebra is isomorphic
 to $\tilde{H}$ defined earlier to Proposition (1.10).  By abuse of 
notation we call the above sub algebra $\tilde{H}$. Now we can not appeal
to Proposition (1.10) as the central operator $t^{\m} K_{k+1}$ has a non
zero degree. But
now by Proposition (1.8) (5), there  exists a proper submodule $W$ for
$\G \otimes A \oplus \Omega_A /d_A \oplus D_k$ such that $V/W$ has finite 
dimensional weight space for $h \oplus \sum \C K_i \oplus D_k$.  Now in
this new module the central operator $t^{\m} K_i$ has zero degree and
leaves finite dimensional weight space invariant.  We now claim that
$t^{\m} K_{k+1}$ is non zero on $V /W$.  Suppose it is zero.  Then $t^{\m}
K_{k+1} V \subseteq W$.  But $V=t^{\m} K_{k+1} V \subseteq W$ a
contradiction as
$W$ is proper.  Thus the $\tilde{H}$ module generated by any $t^{\m}
K_{k+1}$ invariant vector in $V/W$ has all the  properties of the 
Proposition 1.10.  So such a module does not exist.  This proves
$k=n-1$.

\paragraph*{(1.12)}  ~ Let $V$ be an irreducible  $\tau$-module with finite
dimensional weight spaces.  After twisting the module by an automorphism
we can assume that the following central operators are zero.
$$
\begin{array}{lll}
t^{\m} K_i =0 && 1 \leq i \leq n-1 \\
t^m K_n =0 && \ {\rm for} \ m_n \neq 0
\end{array}
$$
Consider $\G_{\rm aff} = \G \otimes \C [t_n, t_n^{-1}] \oplus \C K_n$ with
the obivious Lie bracket.  Consider the Lie algebra homomorphism
$$
\begin{array}{lll}
\varphi = \tau &\to& \G_{\rm off} \ \otimes A_{n-1} \oplus D. \\
X \otimes t^{\underline{r}} &\to&  (X \otimes t_n^{r_n}) \otimes 
t^{\underline{r^1}} \\
&&\underline{r}^1 = (r_1, r_2 \cdots r_{n-1}) \\
\end{array}
$$
$$t^{\m} K_i = \begin{cases} 0 \ {\rm for}  \ 1 \leq i \leq n-1 \\
0 \  i=n, m_n \neq 0 \\
K_n \otimes t^{\m}, m_n=0, i= n
\end{cases}
$$
$$d_i \mapsto d_i$$
It is easy to check that  $\varphi$ defines a Lie algebra homomorphism.
>From the above observation it follows that on an irreducible module
 for $\tau$ with finite dimensional  weight spaces, the Ker $\varphi$
vanishes after twisting by an automorphism.

Thus we can assume $V$  is actually a module for $\G_{{\rm aff}}
\otimes A_{n-1} \oplus D$.  Thus we study only such module.

\section*{Section 2.}

\paragraph*{(2.1)}   We will first introduce some notation

Let $J= \G_{{\rm aff}} \otimes A_{n-1}$ where $\G_{\rm aff} \ = \G
\otimes \C [t_n^{\pm 1}]$ and $A_n=\C [t_1^{\pm 1}, \cdots
t_{n-1}^{\pm 1}]$. Let
$D_k $ be the linear span of $d_{k+1}, \cdots d_n$.  Let $h_k=h \oplus \C k_n
\oplus D_k$ so that $h_{k-1} = h_k \oplus \C d_k$ and $h_0 = h \oplus \C
k_n \oplus D$.

We start with an irreducible $J \oplus D$-module $V$ with finite
dimensional  weight spaces with respect to $h_0$.  We further  assume that
there exists non-zero central operators $z_1, \cdots, z_{n-1}$ with
degree $(k_1, 0 \cdots 0), \cdots (0, \cdots 0, k_{n-1},0)$ without 
loss of generality we can assume that each $k_i$ is a positive integer.
Our aim is this section is to prove that a maximal $J \oplus \C d_n$-sub
module  $W_1$ exists such that $V/W_1$ is an irreducible $J \oplus d_n$-
module with finite dimensional weight spaces with respect to $h_{n-1}$.

We will first produce an irreducible quotient of $V$ for the Lie-algebra 
$J \oplus D_1$ with finite dimensional weight space with respect to $h_1$.

\paragraph*{(2.2)} ~  First notice that $\{z_1 v -v \mid v \in V\}$ is 
a $J \oplus D_1$ proper submodule of $V$.  Thus $V$  is reducible  as
$J \oplus D_1$-module.  So let $W$  be a non-zero proper $J \oplus D_1$
submodule of $V$.  Let $\overline{\mu}$ be a weight of $V$ with respect to
$h_0$.  Let $\mu= \overline{\mu} \mid h_1$.  Let $w \in W$ be any 
$\mu-$weight vector.

Let $w= \displaystyle{\sum_{1 \leq i \leq m}} \ w_i, w_i \in V_{\overline{\mu}+
k_i \delta_i}$ where $k_1 < k_2 < \cdots < k_m$.  Note that $m \geq 2$
as $W$ can not contain $h_0$-weight vectors.  If it contains then it will
not be proper.  Define $d(W)= k_{\ell} -k_1$.  A vector $w$ in $W$ is 
called minimal if $d(W)$ is minimal for  all the $\mu$-weight vector.  Now
one
can prove, that $d(W)$ depends on $W$ but not on $w$.  The proof is same
as
given in Lemma (3.8) of [E2].  (In the proof of the Lemma one has to 
assume $k_{n} -k_1+ \ell_1 \leq \ell_m$ to conclude $d (Xw-v) < 
\ell_m -\ell_1$.  The rest of the proofs works as it is ).
As in the proof of Lemma  (3.9) of [E2], we can prove that $V/W$ has 
finite dimensional weight space for $h_1$.  We can further prove, as in 
Proposition (3.11) of [E2].  $V$ has a maximal $J \oplus D_1$-submodule so
that the quotient is irreducible.  Thus we have proved the following:

\paragraph*{(2.3) Proposition}  Let $V$ be as defined in 2.1.  Then there
exists a maximal $J \oplus D_1$ submodule $W_1$ of $V$ such that
$V/ W$ is irreducible and has finite dimensional weight spaces
with respect to $h_1$.

\paragraph*{(2.4)}  Let $V/ W_1$  as above.  Consider $\{z_2 v-v \mid
v \in V/ W_1\}$.  This is a $J \oplus D_2$ non-zero submodule of
$V/ W_1$ and vectors are not $d_2$ invariant.  But $V/ W_1$
contains $d_2$-invariant vector.  Thus $\{Z_2 v-v \mid v \in V/ W_1\}
$ is a non-zero proper submodule of $V/W_1$.  In particular $V/ W_1$
is $J \oplus D_2$ reducible.  Now by repeating the above process.  We
get an irreducible quotient of $V/ W_1$  which $J \oplus D_2$- module with
finite dimensional weight spaces with respect to $h_2$.  By repeating the
process $n-1$ times we have the following theorem.

\paragraph*{(2.5) ~ Theorem}  Let $V$ be as defined in (2.1).  Then there
exists a maximal $J \oplus D_{n-1}$ submodule $W$ of $V$ so that $V/ W$
is an irreducible $J \oplus D_{n-1}$-module and has finite dimensional
weight spaces with respect to $h_{n-1}$.  Further we can choose $W$ such
that $W$ contains $\{z_iv -v \mid v \in V, i=1, \cdots n-1\}$.

\section*{Section 3}

Let $V$ be irreducible $J \oplus D$-module with finite dimensional weight
spaces with respect to $h_0$.  Assume that there exists non-zero central
operators $z_1, z_2, \cdots, z_{n-1}$ of degree $\underline{k}_1 = (k_1, 0 
\cdots 0) \cdots \underline{k_{n-1}} = (0, \cdots (0, \cdots k_{n-1},0)$.
Let $W$ be a maximal $J \oplus D_{n-1}$ submodule of $V$ so that
$V/ W$ is irreducible $J \oplus D_{n-1}$-module with finite 
dimensional weight space with respect to $h_{n-1}$.  Let $V_1= V/W$.

\paragraph*{(3.1) Lemma}  Each $z$ in $Z$ acts as scalar on $V_1$. 
Further each $z$ in $Z$ which acts as non-zero on $V$ acts as non-zero 
scale on $V_1$.

\paragraph*{Proof} First note that for any $z$ that acts non trivially on
$V$ is of degree $(*, \cdots *, 0)$.  Such  a $z$ leaves each finite
dimensional weight space of $V_1$ invariant.  Thus $z$ has a 
eigen vector with eigen value $\ell$.  Since $V_1$ is irreducible $z$ acts
by the same scalar on the full module $V_1$.   
Suppose $z$ is non-zero on $V$.  Then we claim $\ell \neq 0$.  Suppose 
$\ell=0$ then it follows that $z V \subseteq W$.  But $V=zV \subseteq W$ 
which contradicts the fact that $W$ is proper.  So $\ell \neq 0$.  This 
proves the Lemma.

\paragraph*{(3.2)}  Let $V_1$ be as above.  We now define $J \oplus D$ 
module structure on $V_1 \otimes A_{n-1}$.  Recall  $J= \G_{{\rm aff}}
\otimes A_{n-1}$.  Let $\overline{m} = (m_1, \cdots m_{n-1}) \in \Z^{n-1}$.
Let $X ( \overline{m}) = X \otimes t^{\overline{m}}, X \in \G_{{\rm aff}}, 
t^{\overline{m}} = t_1^{m_1} \cdots t_{n-1}^{m_{n-1}} \in A_{n-1}$.
Define a module structure $\pi (\alpha)$ on $V_1 \otimes A_{n-1}$
where $\alpha = (\alpha_1, \cdots \alpha_{n-1} ) \in \C^{n-1}$.
$$
\begin{array}{lll}
\pi (\alpha) \cdot X (\overline{m}) v (\overline{n}) &=& (X (\overline{m})v)
(\overline{m}+ \overline{n}) \\
d_n v (\overline{n}) &=& (d_n v) (\overline{n}) \\
d_i v (\overline{n}) &=& (\alpha_i + n_i) v (\overline{n}), 1 \leq i 
\leq n-1
\end{array}
$$
The aim of this section is to prove that $(V_1 \otimes A_{n-1}, \pi
(\alpha))$
is completely reducible and $V$ is isomorphic to some component of 
$V_1 \otimes A_{n-1}$ for a suitable $\alpha$.  First consider the
$J \oplus D_{n-1}$ map $S$ from $V_1 \otimes A_{n-1} \to V_1$ given by
 $$S (v (\overline{n})) =v.$$

\paragraph*{(3.3) Lemma} Let $W$ be any non-zero $J \oplus D$-submodule
 of $V_1 \otimes A_{n-1}$. \\
(1) ~~ $S (W) = V_1$ \\
(2) If $\lambda$ is a weight of $V_1$ then $\lambda+
\delta_{\overline{m}}$
is a weight of $W$ for some $\overline{m}$.

\paragraph*{Proof}   Since $V_1$ is irreducible, it is sufficient to prove
that $S(W)$ is non-zero.  But $W$ is $h_0$ weight module and hence contains
vectors of the form $v (\overline{m})$.  Then $S (v (\underline{m})=v \neq
0$.

(2) Let $w \in V_1$ be a weight vector of weight $\lambda$.  Then there 
exists $v \in W $ such that $ S(w) =v$.  Now write $w= \sum w_{\overline{m}}$.
Since $w$ is a weight module it follows that $w_{\overline{m}} \in W$.
Thus (2)  follows.

\paragraph*{(3.4)}  Suppose $z_i$  acts by $\ell_i $ on $V_1$ and we know that
$\ell_i \neq 0$.

Consider $z_i v (\overline{m}) = \ell_i v (\overline{m} +\underline{k}_i)$.
Clearly  $z_i$ is central operator and invertible on $V_1 \otimes A_{n-1}$.  
Let 
$$\Gamma =k_1 \Z \oplus  k_2 \Z \oplus  \cdots \oplus k_{n-1} \Z 
\subseteq  \Z^{n-1}.$$  
It is easy to see that for each $\overline{m} \in
\Gamma$ there exists a  central operator $z_{\overline{m}}$ on $V_1
\otimes
A_{n-1}$ which is invertible.

Let $W$ be a submodule of $V_1 \otimes A_{n-1}$.  

We claim then dim $W_{\lambda}=$ dim $W_{\lambda}+ \delta_{\overline{m}}$
for
all $\underline{m} \in \Gamma$.  Note that $z_{\overline{m}} W_{\lambda} 
\subseteq W_{\lambda+ \delta_{\overline{m}}}$ and equality holds as 
$z_{\overline{m}}$ is invertible.  Thus the claim follows.  let $\lambda$
be a fixed weight of $V_1$.  Let $K(W)= 
\displaystyle{\sum_{0 \leq m_i < k_i}} \  $ dim $ (W_{\lambda + 
\delta_{\underline{m}}})$.

Now from Lemma 3.3 (2) and the above claim it follows that $K(W) \neq 0$.

\paragraph*{(3.5) Proposition}  $V_1 \otimes A_{n-1}$ contains an 
irreducible submodule.

\paragraph*{Proof}  We will prove that for any non-zero submodule $W$ of
$V_1 \otimes A_{n-1}$ either $W$ contains an irreducible submodule or
$W$ contains a  module $W_1 \ni K (W_1) <K (W)$.  Fix a non-zero
submodule $W$ of $V_1 \otimes A_{n-1}$.  Consider $S= 
\displaystyle{\oplus_{\overline{m} \in \Z^{n-1}}} \  W_{\lambda+
\delta_{\overline{m}}}$.  Let $\tilde{W}$ be a submodule of $W$ generated 
by $S$.  Then clearly $K (\tilde{W}) = K (W)$.  Suppose $\tilde{W}$ has
a proper submodule $W_1$.

\paragraph*{Claim}  $K(W_1) <K (\tilde{W})$.  Suppose the claim is
false.  Consider
$$W_1 \cap (V_1 \otimes A_{n-1})_{\lambda+ \delta_{\overline{m}}} 
\subseteq \tilde{W} \cap (V_1 \otimes A_{n-1})_{ \lambda+
\delta_{\overline{m}}}. 
\leqno{(*)}$$  
Now $K (W_1) = K (\tilde{W})$ implies that equality holds for all $\overline{m}
\ni 0 \leq m_i < k_i$.  Now by the claim in (3.4) it follows that 
equally holds in * for all $\overline{m}$ in $ \Z^{n-1}$.  This proves
$W_1$
contain
$S$ and hence contain $\tilde{W}$.  This contradicts the fact that $W_1$
is proper in $\tilde{W}$.  Hence the claim.  Now by repeating the above
process we get an irreducible submodule of $V_1 \otimes A_{n-1}$.  We now 
prove complete reducibility of $V_1 \otimes A_{n-1}$.

\paragraph*{(3.6)} ~ Let $W$ be an irreducible submodule of $V_1 \otimes
A_{n-1}$.  We have fixed a weight $\lambda$ of $V_1$.  By Lemma 3.3 (b)
there exists  a weight vector $v$ in $V$ of weight $\lambda$ such that
$v (\overline{r}) \in W$.

Let $\overline{U}$ be the universal enveloping algebra of $J \oplus D$.
Write
$$\overline{U}= \displaystyle{\otimes_{\overline{m} \in \Z^{n-1}}} \
\overline{U}_{\overline{m}}$$
$$\overline{U}_{\overline{m}}=\{X \in U \mid
[d_i,X]=r_i X, 1 \leq i \leq n-1\}$$
Clearly $\overline{U} v (\overline{r}) =W$ where
$\overline{U} v (\overline{r})$ is the $J \oplus D$-module 
generated by $v (\overline{r})$

\paragraph*{(3.7) ~ Lemma} (1) $\overline{U} v (\overline{s})$ is
irreducible $J \oplus D$-module for any $\overline{s} \in \Z^{n-1}$.

(2) $\displaystyle{\sum_{0 \leq s_i<k_i}} \ \overline{U} v (\overline{s})
= V_1 \otimes A_{n-1}$.

\paragraph*{Proof}  Consider the map $\varphi: \overline{U} v 
(\overline{r}) \to \overline{U}v (\overline{s})$ given by
$\varphi (w (\underline{k}))= w (\overline{k}+ \overline{s}-
\overline{r})$.
It is a $J \oplus D_{n-1}$ module map.  But need not be $J \oplus D$-
module map.  Nevertheless it is an isomorphism.  It is easy to see that
$\overline{U} v (\overline{s})$ is irreducible.

\paragraph*{Claim} $\displaystyle{\sum_{\overline{s} \in \Z^{n-1}}}
\overline{U} v (\overline{s}) =V_1
\otimes A_{n-1}$.  Let $w (\underline{s}) \in V_1 \otimes A_{n-1}, w \in V_1$.
Since $V_1$ is irreducible there exists $X \in \overline{U} \ni X v = w$.
Write
$X = \sum X_{\overline{r}}$ where $X_{\overline{r}} \in 
\overline{U}_{\overline{r}}$.

Consider
$$\sum X_{\overline{r}} \cdot v (\overline{s}- \overline{r}) =\sum
(X_{\overline{r}} v) (\overline{s}) = w (\overline{s}).$$
This proves the claim.  Now recall that there exists a non-zero central
operator $z_{\overline{m}}$ for any $\underline{m} \in \Gamma$.  Further
$z_{\overline{m}}$ is invertible.

Thus $z_{\overline{m}} \overline{U} v (\overline{r}) = U (\overline{r}+
\overline{m})$.  From this it follows that
$$\displaystyle{\sum_{\overline{s} \in \Z^{n-1}}} \overline{U} v 
(\overline{s}) = \displaystyle{\sum_{0 \leq s_i < k_i}}  \overline{U} v
(\overline{s}).$$

\paragraph*{(3.8)  Proposition}  $V_1 \otimes A_{n-1}$ is completely 
reducible as $J \oplus D$-module.

\paragraph*{Proof} We have already  seen that $\displaystyle{\sum_{0 \leq s_i
< k_i}} \overline{U}  v (\overline{s}) = V_1 \otimes A_{n-1}$.
The sum is finite.  Let $T= \{ \underline{s} \mid 0 \leq s_i< k_i \}$.  
Suppose $\overline{U} v (\overline{s}) \cap 
\displaystyle{\sum_{\stackrel{\underline{r} \in T}{\underline{r} \neq 
\underline{s}}}} \overline{U} v (\overline{s}) \neq 0$.
Since $\overline{U} v (\overline{s})$ is irreducible it follows that
$\overline{U} v (\overline{s}) \subseteq \displaystyle{\sum_{\stackrel{
\overline{r} \in T}{\overline{s} \neq \overline{r}}}} \overline{U} v 
(\overline{r})$. Thus we have
reduced the finite sum by 1.  Repeating this process we can replace the 
$\sum$ by direct sum with fewer terms.

Let $\overline{U}= \displaystyle{\oplus_{\eta \in h_0^*}} \overline{U}_{\eta},
\overline{U}_\eta = \{X \in \overline{U} \mid [h, X] = \eta (h) X \forall
h \in h_0 \}$.  Let $U$ be the universal enveloping algebra of $J \oplus
D_n$.
Let ${U}= \displaystyle{\oplus_{\eta \in h_{\eta \in h_{n-1}^{*}}}}
U_{\eta}, U_{\eta} = \{ X \in U \mid [h, X]= \eta (h) X \  
\forall \in h_{n_1}\}$
\paragraph*{(3.9) Theorem } $V$ is isomorphic to a component of $V_1 \otimes
A_{n-1}$ as $J \oplus D$-module.

\paragraph*{Proof}  It can be proved similar to the techniques of 3.13, 3.14,
3.15 and  3.16 of [E2].  The $U$ and $\overline{U}$ is taken as defined by
us. The $d$ in [E2] is taken as the space spanned by $d_1, \cdots
d_{n-1}$.  We do not need the condition (2.2) of [E2].  In stead we have
non-zero central operators $z_1, \cdots, z_{n-1}$ which are invertible in
$V$ as well as $V_1 \otimes A_{n-1}$.

\paragraph*{(3.10) Remark}  The study of irreducible $\tau$-modules with
finite dimensional weight space where some part of the center acts
non-trivially is now reduced to the study of irreducible modules of 
$J \oplus D_n$ with finite dimensional weight spaces.  The advantage of
the later modules is that the infinite dimensional center acts as
scalars.

\begin{center}
{\bf REFERENCES}
\end{center}
\vskip 5mm

\begin{enumerate}
\item[{[BC]}] Berman, S. and Cox, B. Enveloping algebras and representations
of toroidal Lie-algebras, Pacific Journal of Math. 165 (1994), 239-267.
\item[{[CO]}] Cox, Ben. Two Realizations of Toroidal $s \ell_2 (\C)$,
Contemporary Mathematics, 297, (2002), 47-68.
\item[{[CT]}] Chari, V. and Thang Le, Representations of double affine
Lie-
algebras, Preprint, arxiv / org - 0205312.
\item[{[E1]}] Eswara Rao, S. Iterated loop modules and  a filteration for
vertex representation of toroidal Lie-algebras, Pacific Journal of Mathematics,
171 (2), (1995), 511-528.

\item[{[E2]}] Eswara Rao, S. Classification of loop modules with finite
dimensional weight spaces, Math. Anal., 305 (1996), 651-663.
\item[{[E3]}] Eswara Rao, S. Classification of irreducible integrable
modules for multi-loop algebras with finite dimensional weight spaces,
Journal of Algebra, 246, (2001), 215-225.
\item[{[E4]}] Eswara Rao, S. Classification of irreducible integrable
modules for toroidal Lie-algebras with finite dimensional weight spaces,
TIFR Preprint (2001) arxiv / org. 0209060.
\item[{[EM]}] Eswara Rao, S. and Moody, R.V. Vertex representations for $n$-
toroidal Lie-algebras and a generalization of Virasoro algebra, Comm.
Math. Physics, 159 (1994), 239-264.
\item[{[EMY]}] Eswara Rao, S., Moody, R.V., and Yokonuma, T. Toroidal
Lie-algebras and vertex representations, Geom. Ded. 35, (1990), 283-307.
\item[{[FK]}] Futorny, V. and Kashuba, I. Verma type modules for 
toroidal Lie algebras, Comm.in Math. 27 (8), (1999), 3979-3991.
\item[{[F]}] Futorny, V. Representations of affine Lie-algebras, Queen's 
papers in Pure and Applied Mathematics, (1997).
\item[{[JMT]}] Jing, N., Misra, K. and Tan, S. Bosonic realization of
higher level toroidal Lie-algebra, Preprint 2002.
\item[{[K]}] Kassel, C., Kahler differentials and coverings of complex
simple Lie-algebra extended over a commutative algebra, J Pure Appl.
Algebra, 34 (1985), 266-275.
\item[{[MY]}] Morita, J. and Yoshi, Y. Universal central extensions of
Chevally Algebras over Laurent polynomials rings and GIM Lie algebra, Proc.
Japan Acad. Ser.A 61 (1985), 179-181.
\item[{[Y]}] Youngsun, Yoon. On the polynomial representations of 
current algebras, Journal of Algebra, 252 (2) (2002), 376-393.

\end{enumerate}

\end{document}